\documentclass[reqno]{amsart}
\usepackage{amssymb,latexsym,amsmath,amsthm,enumerate,amsbsy}
\usepackage[mathscr]{eucal}
\usepackage{color,graphicx}
\usepackage{mathrsfs}
\usepackage{tikz}
\usepackage{cite}
\usepackage{longtable}
\usepackage{caption} 
\usepackage{etoolbox} 
\usepackage{tabularx}
\usepackage{array}

\usepackage{enumitem}

\usepackage{booktabs} 

\makeatletter
\@namedef{subjclassname@2020}{\textup{2020} Mathematics Subject Classification}
\makeatother

\usetikzlibrary{positioning,decorations.pathreplacing,patterns,decorations.pathmorphing}
\tikzset{%
element/.style={draw, shape=circle, fill=white, inner sep=1.4pt}
}

\theoremstyle{plain}
\newtheorem{thm}{Theorem}[section]

\newtheorem{prob}[thm]{Problem}

\theoremstyle{definition}

\newcommand{\M}{\mathbf M}

\newcounter{eqbasis}
\newcommand{\bq}{\mathbf{q}}
\newcommand{\br}{\mathbf{r}}
\newcommand{\bs}{\mathbf{s}}

\newcommand{\bu}{\mathbf{u}}
\newcommand{\bv}{\mathbf{v}}
\newcommand{\bw}{\mathbf{w}}

\theoremstyle{plain}
\newtheorem{theorem}[thm]{Theorem}
\newtheorem{lemma}[thm]{Lemma}
\newtheorem{proposition}[thm]{Proposition}
\newtheorem{corollary}[thm]{Corollary}
\theoremstyle{definition}

\newcommand{\V}{\mathsf V}

\begin{document}
\title[A nonfinitely based ai-semiring with finitely based matrix semirings]
{A nonfinitely based additively idempotent semiring with finitely based matrix semirings}
\author{Jun Jiao}
\address{School of Mathematics, Northwest University, Xi'an, 710127, Shaanxi, P.R. China}
\email{jjunjiao@163.com}

\author{Miaomiao Ren}
\address{School of Mathematics, Northwest University, Xi'an, 710127, Shaanxi, P.R. China}
\email{miaomiaoren@yeah.net}



\begin{abstract}
We establish that the matrix semiring $\mathbf{M}_n(SR_6)$
over the six-element nonfinitely based additively idempotent semiring $SR_6$
is finitely based for every $n \geq 2$.
This provides the first known example of an additively idempotent semiring
$S$ that is nonfinitely based, while the matrix semirings $\mathbf{M}_n(S)$
are finitely based.
Moreover, we show that the interval $[\mathsf{V}(SR_6), \mathsf{V}(\mathbf{M}_2(SR_6))]$
in the lattice of additively idempotent semiring varieties has the cardinality of the continuum.
\end{abstract}

\maketitle

\section{Introduction}\label{sec:intro}
An \emph{additively idempotent semiring} (or \emph{ai-semiring} for short)
is an algebra $(S, +, \cdot)$ such that the additive reduct $(S, +)$
is a commutative idempotent semigroup, the multiplicative reduct
$(S, \cdot)$ is a semigroup, and the distributive laws
\[
x(y+z)\approx xy+xz,\quad (x+y)z\approx xz+yz
\]
hold.
Well-known examples of ai-semirings include max-plus algebras~\cite{aei},
flat extensions of groups~\cite{jackson08}, and distributive lattices~\cite{burris1981}.
These and other similar algebras provide a rich algebraic framework for many areas,
forming the foundation of tropical geometry~\cite{ms} and tropical linear algebra~\cite{but},
with important applications in theoretical computer science~\cite{go} and information science~\cite{gl}.

Let $S$ be an ai-semiring. Then $(S, +)$ is a semilattice.
The binary relation $\leq_S$ on $S$ defined by
\[
a \leq_S b \Leftrightarrow a+b=b
\]
is a partial order on $S$, often referred to simply as $\leq$.
One readily checks that $\leq$ is compatible with both addition and multiplication.
Moreover, the supremum of $a$ and $b$ with respect to $\leq$ is $a+b$.
In particular, the addition of $S$ is determined by $\leq$.
For this reason, an ai-semiring is often called a \emph{semilattice-ordered semigroup}.
Consequently, when describing an ai-semiring, it is often convenient to specify $\leq$
by means of its Hasse diagram rather than to give the additive Cayley table.
Whenever an order on an ai-semiring is mentioned, it always refers to the order defined above.

A \emph{variety} of ai-semirings is a class of ai-semirings that is closed
under taking subalgebras, homomorphic images, and arbitrary direct products.
By Birkhoff's theorem, a class of ai-semirings is a variety if and only if
it is an equational class; that is, it is the class of all ai-semirings
satisfying a certain set of identities. Such a set is called an
\emph{equational basis} of the variety.

An ai-semiring variety is \emph{finitely based} if it admits a finite
equational basis; otherwise, it is \emph{nonfinitely based}.
An ai-semiring $S$ is \emph{finitely based} (\emph{nonfinitely based})
if the variety $\V(S)$ generated by $S$ is finitely based (nonfinitely based).
If $\Sigma$ is an equational basis of $\V(S)$,
we also say that $\Sigma$ is an equational basis of $S$.

The finite basis problem for a class of ai-semirings, one of the most important
topics in universal algebra, concerns the classification of its members
with respect to the finite basis property.
Over the past two decades, the finite basis problem for ai-semirings has received considerable attention; 
see, for example,~\cite{aei, gr2026, Dolinka2009, jackson08, jrz, jr2025, JiaoRen2026b, LRY, ShaoRen2015, Volkov}.

Let $n\geq 2$ be an integer, let $S$ be an ai-semiring, and let
$\mathbf{M}_n(S)$ denote the set of all $n\times n$ matrices over $S$.
Then $\M_n(S)$ forms an ai-semiring under the usual matrix addition and multiplication; that is,
\[
A+B=[a_{ij}+b_{ij}]_{n\times n},\quad
AB=\left[\sum_{k=1}^n a_{ik}b_{kj}\right]_{n\times n}
\]
for all $A=[a_{ij}]_{n\times n}, B=[b_{ij}]_{n\times n}\in\M_n(S)$.
This ai-semiring is called the \emph{matrix semiring over $S$}.
It is straightforward to verify that $A \leq B$ if and only if $a_{ij}\leq b_{ij}$ for all $1 \leq i, j \leq n$.

As observed in \cite{jr2025, JiaoRen2026b}, one can obtain the following infinite strict ascending chain
\[
S\hookrightarrow \mathbf{M}_2(S)\hookrightarrow \mathbf{M}_{3}(S)
\hookrightarrow\cdots \hookrightarrow \mathbf{M}_n(S)\hookrightarrow \mathbf{M}_{n+1}(S)\hookrightarrow \cdots,
\]
and hence the following infinite ascending chain of varieties:
\[
\V(S) \leq \V(\mathbf{M}_2(S))  \leq  \V(\mathbf{M}_{3}(S))
 \leq \cdots  \leq   \V(\mathbf{M}_n(S))  \leq  \V(\mathbf{M}_{n+1}(S))  \leq  \cdots.
\]

Despite the considerable progress in the study of the finite basis problem
for ai-semirings, the corresponding problem for matrix semirings over
ai-semirings remains largely unknown.
Some progress has been made in this direction; see, for example, \cite{Dolinka2009, jr2025, JiaoRen2026b, Volkov}.
To the best of our knowledge,
Dolinka~\cite{Dolinka2009} was the first to consider the finite basis problem
for matrix semirings over ai-semirings.
More precisely, he showed that every locally finite variety containing
$\mathbf{M}_2(D_2)$ is nonfinitely based, where $D_2$
denotes the two-element distributive lattice.
Consequently, $\mathbf{M}_n(D_2)$ is nonfinitely based for all $n \geq 2$.
Since $D_2$ is finitely based (see \cite{ShaoRen2015}), this shows that
there exists a finitely based ai-semiring $S$ such that the matrix semiring
$\mathbf{M}_n(S)$ is nonfinitely based for some $n\geq 2$.

Jiao and Ren~\cite{jr2025} provided a complete classification of matrix semirings
over two-element ai-semirings with respect to the finite basis property.
They proved that for every integer $n \geq 2$ and every two-element ai-semiring $S$,
$\mathbf{M}_n(S)$ is finitely based if and only if $S$ is distinct from $D_2$.
Since every two-element ai-semiring $S$ is finitely based (see \cite{ShaoRen2015}),
it follows that the matrix semirings $\mathbf{M}_n(S)$ preserve the finite basis property
whenever $S$ is a two-element ai-semiring distinct from $D_2$.

Motivated by these results, at the end of \cite{jr2025}, Jiao and Ren posed the following natural problem.

\begin{prob}\label{prob:matrix}
Does there exist a nonfinitely based ai-semiring $S$ and an integer
$n \geq 2$ such that $\mathbf{M}_n(S)$ is finitely based?
\end{prob}

In contrast, Jiao and Ren~\cite{JiaoRen2026b} proved that the matrix semiring
$\mathbf{M}_n(S_7)$ is nonfinitely based for all $n \geq 2$, where $S_7$ is
the particular three-element ai-semiring defined in \cite{jrz}.
Since $S_7$ itself is nonfinitely based (see \cite{jrz}), this confirms that
the nonfinite basis property is indeed preserved by matrix constructions in this case.

At the time, it seemed highly counterintuitive that there should exist a nonfinitely based ai-semiring
$S$ such that $\mathbf{M}_n(S)$ is finitely based;
one would naturally expect matrix constructions to preserve the nonfinite basis property.
However, the main result of this paper provides an affirmative answer to Problem~\ref{prob:matrix}.

We now present a concrete example that achieves this.
Let $SR_6$ denote the six-element ai-semiring $\{\infty, \omega, a, b, c, d\}$ whose Hasse diagram
and multiplicative Cayley table
are presented in Figure~\ref{fig:sr6} and Table~\ref{tab:sr6}, respectively.
This algebra possesses a remarkable property.
More precisely, Lyu, Ren and Yue~\cite[Theorem 5.4]{LRY} showed that
the variety $\mathsf{V}(SR_6)$ itself is nonfinitely based,
while all proper subvarieties of $\mathsf{V}(SR_6)$ are finitely based.
In fact, the subvariety lattice of $\mathsf{V}(SR_6)$ is a four-element chain \cite[Proposition 5.3]{LRY}.
Consequently, $SR_6$ is nonfinitely based.

\begin{figure}[ht]
\refstepcounter{figure}
\label{fig:sr6}
\noindent
\parbox{0.48\textwidth}{\small Fig. 1 The additive order of $SR_6$}
\hfill
\begin{tikzpicture}[
  baseline=(inf-node.base),
  x=1.0cm,y=0.7cm,
  element/.style={font=\normalsize,inner sep=2pt,fill=white}
]
\coordinate (inf) at (0,4);
\coordinate (d)   at (-1.5,3);
\coordinate (a)   at (1.5,3);
\coordinate (b)   at (-1.5,2);
\coordinate (om)  at (0,2);
\coordinate (c)   at (1.5,2);

\draw (b) -- (d) -- (inf);
\draw (om) -- (inf);
\draw (c) -- (a) -- (inf);

\node[element] at (inf) (inf-node) {$\infty$};
\node[element] at (d)   {$a$};
\node[element] at (a)   {$d$};
\node[element] at (b)   {$c$};
\node[element] at (om)  {$\omega$};
\node[element] at (c)   {$b$};
\end{tikzpicture}
\end{figure}

\noindent
\raisebox{-1.5ex}{\parbox[t]{0.48\textwidth}{\small Tab. 1 The multiplication table of $SR_6$}}
\refstepcounter{table}
\label{tab:sr6}
\hfill
{
\setlength{\tabcolsep}{8pt}
\renewcommand{\arraystretch}{1.1}
\begin{tabular}[t]{c c c c c c c}
\toprule
$\cdot$ & $\infty$ & $\omega$ & $a$ & $b$ & $c$ & $d$ \\
\midrule
$\infty$ & $\infty$ & $\infty$ & $\infty$ & $\infty$ & $\infty$ & $\infty$ \\
$\omega$ & $\infty$ & $\infty$ & $\infty$ & $\infty$ & $\infty$ & $\infty$ \\
$a$      & $\infty$ & $\infty$ & $\infty$ & $\omega$ & $\infty$ & $\infty$ \\
$b$      & $\infty$ & $\infty$ & $\omega$ & $\infty$ & $\omega$ & $\infty$ \\
$c$      & $\infty$ & $\infty$ & $\infty$ & $\omega$ & $\infty$ & $\omega$ \\
$d$      & $\infty$ & $\infty$ & $\infty$ & $\infty$ & $\omega$ & $\infty$ \\
\bottomrule
\end{tabular}
}

\vspace{12pt}

Note that the underlying set of $SR_6$ in \cite{LRY}
is labelled differently from ours; the two algebras are isomorphic,
and we identify them throughout.
We are now ready to state our main result.

\begin{theorem}\label{maintheorem}
Let $n \geq 2$ be an integer. Then the matrix semiring $\mathbf{M}_n(SR_6)$ is finitely based.
\end{theorem}

Theorem~\ref{maintheorem} follows directly from Proposition~\ref{sr:thm:matrix},
which provides an explicit finite equational basis for $\mathbf{M}_n(SR_6)$.
Since $SR_6$ is nonfinitely based,
Theorem~\ref{maintheorem} resolves Problem~\ref{prob:matrix} and provides
the first example of a nonfinitely based ai-semiring $S$ for which
$\mathbf{M}_n(S)$ is finitely based, even for all $n \geq 2$.
The paper is organized as follows. In Section~\ref{sec:preli}, we collect the necessary preliminaries.
In Section 3, we present a finite equational basis for $\mathbf{M}_n(SR_6)$
and show that the interval $[\mathsf{V}(SR_6), \mathsf{V}(\mathbf{M}_2(SR_6))]$
has the cardinality of the continuum.

\section{Preliminaries}\label{sec:preli}
In this section, we recall some basic terminology and notation that will be used throughout the paper.

Let $X$ be a set of variables, and let $X^+$ be the free semigroup over $X$.
An \emph{ai-semiring term} (or simply a \emph{term}) over $X$ is
a formal sum of words in $X^+$, where the order of the summands is irrelevant
and repeated occurrences of the same word are identified with a single occurrence.
In what follows, terms and words are denoted by bold lowercase letters $\mathbf{u}, \mathbf{v}, \mathbf{w}, \dots$,
while ordinary lowercase letters $x, y, z, \dots$ stand for variables.
Consequently, every term can be written as $\bu_1 + \bu_2 + \cdots + \bu_n$,
where $\bu_1, \bu_2, \ldots, \bu_n$ are distinct words in $X^+$.
For convenience, we sometimes identify a term with its underlying set of words;
thus, we write $\bw \in \bu$ to mean that the word $\bw$ is a summand of the term $\bu$.

An \emph{ai-semiring identity} (or simply an \emph{identity}) is a
formal expression of the form $\bu\approx \bv$, where $\bu$ and $\bv$ are terms.
Let $S$ be an ai-semiring, and let $\bu\approx \bv$ be an identity over
$\{x_1, x_2, \ldots, x_n\}$.
We say that $S$ \emph{satisfies} $\bu\approx \bv$, or that
$\bu \approx \bv$ \emph{holds} in $S$, if
\[
\bu(a_1, a_2, \ldots, a_n)=\bv(a_1, a_2, \ldots, a_n)
\]
for all $a_1, a_2, \ldots, a_n \in S$. Here, $\bu(a_1, a_2, \ldots, a_n)$ denotes
the value of $\bu$ in $S$ under the mapping $x_i\mapsto a_i$ for $1\leq i \leq n$,
and similarly for $\bv(a_1, a_2, \ldots, a_n)$.
Such a mapping is called an \emph{assignment into $S$}.
Therefore, $S$ satisfies $\bu\approx \bv$ if and only if $\varphi(\bu)=\varphi(\bv)$
for every assignment $\varphi$ into $S$, where $\varphi(\bu)$ and $\varphi(\bv)$
denote the values of $\bu$ and $\bv$ in $S$ under $\varphi$, respectively.

We write $\bu \preceq \bv$ (or equivalently, $\bv \succeq \bu$)
to denote the identity $\bv \approx \bv+\bu$, and call it an
\emph{ai-semiring inequality} (or simply an \emph{inequality}).
Let $\bu\approx \bv$ be an identity such that
\[
\bu=\bu_1+\cdots+\bu_k,\quad \bv=\bv_1+\cdots+\bv_\ell,
\]
where $\bu_i, \bv_j\in X^+$ for $1\leq i\leq k$ and $1\leq j\leq \ell$.
It is readily verified that the ai-semiring variety defined by $\bu\approx \bv$
is exactly the same as the variety defined by the inequalities
\[
\bu_i\preceq \bv, \quad \bv_j\preceq \bu \qquad (1\leq i\leq k, \ 1\leq j\leq \ell).
\]
Consequently, throughout this paper, we restrict our attention to inequalities
of the form $\bq\preceq\bu$, where $\bq$ is a word and $\bu$ is a term.

To end this section, we introduce some notation.
For a nonempty word $\bw$, let $c(\bw)$ denote its \emph{content}; that is,
the set of all variables occurring in $\bw$, and let $\ell(\bw)$ denote its
\emph{length}; that is, the number of variables occurring in $\bw$ counted with multiplicity.

Now let $\bu$ be a term such that $\bu=\bu_1+\bu_2+\cdots+\bu_n$,
where $\bu_i \in X^+$ for $1 \leq i \leq n$.
Then $c(\bu)$ denotes the \emph{content} of $\bu$; that is,
\[
c(\bu)=\bigcup_{i=1}^n c(\bu_i).
\]
For a positive integer $k$, let $L_{k}(\bu)$ denote the term consisting of
the sum of all $\bu_i$ such that $\ell(\bu_i)=k$.

\section{Finite equational basis for $\M_n(SR_6)$}\label{sec:finitecounter}
In this section, we first present an explicit finite equational basis for
$\mathbf{M}_n(SR_6)$ for $n \geq 2$, which immediately implies that
every $\mathbf{M}_n(SR_6)$ is finitely based.
We also show that the interval $[\mathsf{V}(SR_6), \mathsf{V}(\mathbf{M}_2(SR_6))]$
has the cardinality of the continuum.
As a first step, we collect some basic properties about $SR_6$.

\begin{lemma}\label{lemma26111501}
The following statements about $SR_6$ hold:
\begin{enumerate}[label=$(\alph*)$]
\item \label{item1}
The element $\infty$ is the additive maximum element of $SR_6$;

\item \label{item2}
The element $\infty$ is
the multiplicative zero element of $SR_6$;

\item \label{item4}
The product of any two elements of $SR_6$ is either $\infty$ or $\omega$;

\item \label{item3}
The product of any three elements of $SR_6$ equals $\infty$;

\item \label{item5}
$p+\omega=\infty$ for all $p \in SR_6\setminus \{\omega\}$;

\item \label{item6}
$p\omega=\omega p=\infty$ for all $p \in SR_6$.
\end{enumerate}
\end{lemma}
\begin{proof}
This follows directly from Figure~\ref{fig:sr6} and Table~\ref{tab:sr6}.
\end{proof}

As a consequence of the properties listed in Lemma~\ref{lemma26111501},
we present three identities of $\M_n(SR_6)$.

\begin{proposition}\label{sr:lem:sound}
Let $n\geq2$ be an integer. Then $\M_n(SR_6)$ satisfies the identities
\begin{align}
&x \preceq x_1x_2x_3; \label{sr:eq:N}\\
& x+xy\approx x_1x_2x_3; \label{sr:eq:L}\\
& x+yx\approx x_1x_2x_3. \label{sr:eq:R}
\end{align}
\end{proposition}
\begin{proof}
Let $A=[a_{ij}], B=[b_{ij}]$, $C=[c_{ij}]$, and $D=[d_{ij}]$ be arbitrary matrices in $\M_n(SR_6)$.
For all $1\leq i, j\leq n$,
\[
(ABC)_{ij}=\sum_{k=1}^n\sum_{l=1}^n a_{ik}b_{kl}c_{lj}=\infty,
\]
which follows from Lemma~\ref{lemma26111501}\ref{item1} and \ref{item3}.
This shows that $ABC$ is the constant matrix $[\infty]_{n\times n}$ with all entries equal to $\infty$.
By Lemma~\ref{lemma26111501}\ref{item1} again, we have $D\leq [\infty]_{n\times n}$,
and so $D\leq ABC$. Consequently, the inequality \eqref{sr:eq:N} holds in $\M_n(SR_6)$.

To show that the identity \eqref{sr:eq:L} holds in $\M_n(SR_6)$,
by the preceding argument, we need to prove that $(A+AB)_{ij}=\infty$.
Indeed, if $A_{ij}=\infty$, then
\[
(A+AB)_{ij}=A_{ij}+(AB)_{ij}=\infty+(AB)_{ij}=\infty,
\]
where the last equality follows from Lemma~\ref{lemma26111501}\ref{item1}.
If $A_{ij}=\omega$, then by Lemma~\ref{lemma26111501}\ref{item1} and \ref{item6},
$(AB)_{ij}=\infty$. Hence
\[
(A+AB)_{ij}=A_{ij}+(AB)_{ij}=A_{ij}+\infty=\infty,
\]
where the last equality follows from Lemma~\ref{lemma26111501}\ref{item1}.
If $A_{ij} \in \{a, b, c, d\}$, then by Lemma~\ref{lemma26111501}\ref{item4}, $(AB)_{ij} \in \{\infty, \omega\}$.
Thus
\[
(A+AB)_{ij}=A_{ij}+(AB)_{ij}=\infty,
\]
where the last equality follows from Lemma~\ref{lemma26111501}\ref{item1} and \ref{item5}.

Similarly, the identity \eqref{sr:eq:R} also holds in $\M_n(SR_6)$.
\end{proof}

We now proceed to give a complete description of the identities of $\M_n(SR_6)$.

\begin{proposition}\label{sr:lem:absorption}
Let $n \geq 2$ be an integer, and
let $\bq \preceq \bu$ be a nontrivial inequality such that
$\bu=\bu_1+\bu_2+\cdots+\bu_m$ with $\bu_i, \bq \in X^+$ for $1\leq i \leq m$.
Then $\bq\preceq \bu$ is satisfied by $\M_n(SR_6)$ if and only if
either $\ell(\bu_i)\geq3$ for some $\bu_i\in \bu$, or $c(L_1(\bu))\cap c(L_2(\bu))\neq \emptyset$.
\end{proposition}
\begin{proof}
Suppose that $\ell(\bu_i) \geq 3$ for some $\bu_i\in \bu$. By applying the inequality \eqref{sr:eq:N}, we obtain
\[
\bu \succeq \bu_i \succeq \bq.
\]
Hence $\bq\preceq \bu$ is derivable from \eqref{sr:eq:N}.
Since \eqref{sr:eq:N} holds in $\mathbf{M}_n(SR_6)$ by
Proposition~\ref{sr:lem:sound}, the inequality $\bq\preceq \bu$ is also
satisfied by $\mathbf{M}_n(SR_6)$.

Assume that $c(L_1(\bu))\cap c(L_2(\bu))\neq \emptyset$.
Then there exist variables $x$ and $y$ such that
either $x, xy \in \bu$ or $x, yx \in \bu$.
If $x, xy \in \bu$, then
\[
\bu \succeq x+xy \stackrel{{\eqref{sr:eq:L}}}\approx x_1x_2x_3 \stackrel{{\eqref{sr:eq:N}}} \succeq \bq.
\]
Therefore, $\bq\preceq \bu$ is derivable from \eqref{sr:eq:L} and \eqref{sr:eq:N}.
Similarly, if $x, yx \in \bu$, then $\bq\preceq \bu$ is derivable from \eqref{sr:eq:R} and \eqref{sr:eq:N}.
Since \eqref{sr:eq:N}, \eqref{sr:eq:L} and \eqref{sr:eq:R} hold in
$\mathbf{M}_n(SR_6)$ by Proposition~\ref{sr:lem:sound}, we conclude that
$\bq\preceq \bu$ is satisfied by $\mathbf{M}_n(SR_6)$.

Conversely, suppose that $\ell(\bu_i)\leq 2$ for all $\bu_i\in \bu$ and $c(L_1(\bu))\cap c(L_2(\bu))=\emptyset$.
To show that $\bq \preceq \bu$ fails in $\mathbf{M}_n(SR_6)$,
it suffices to prove that $\eta(\bq) \nleq \eta(\bu)$ for some assignment $\eta$ into $\M_n(SR_6)$.

For $\alpha, \beta\in SR_6$, define the matrix
\[
M(\alpha, \beta)=
\begin{bmatrix}
\alpha & \beta & \cdots & \beta\\
c      & c     & \cdots & c\\
\vdots & \vdots&        & \vdots\\
c      & c     & \cdots & c
\end{bmatrix}_{n\times n}.
\]
For all $\alpha, \beta, \gamma, \delta\in SR_6$, we have
$M(\alpha, \beta)_{12}=\beta$ and
\begin{equation}\label{id26111901}
\bigl(M(\alpha, \beta)M(\gamma, \delta)\bigr)_{12}
=\alpha\delta+\sum_{k=2}^{n}\beta c
=\alpha\delta+\beta c,
\end{equation}
where the second equality follows from additive idempotence.

\textbf{Case 1.} $\ell(\bq)=1$ or $\ell(\bq)\geq3$.
Define an assignment $\psi$ into $\M_n(SR_6)$ by
\[
\psi(z)=
\begin{cases}
M(c, \omega), &z\in c(L_1(\bu)),\\
M(c, b), &z\notin c(L_1(\bu)).
\end{cases}
\]
If $z\in c(L_1(\bu))$, then $\psi(z)_{12}=\omega$.
If $st\in L_2(\bu)$, then $s, t\notin c(L_1(\bu))$, and hence
\[
\psi(st)_{12}=(\psi(s)\psi(t))_{12}=(M(c,b)M(c,b))_{12}\stackrel{\eqref{id26111901}}=cb+bc=\omega+\omega=\omega.
\]
It follows that
\begin{align*}
\psi(\bu)_{12} &= \psi(L_1(\bu)+L_2(\bu))_{12}
= (\psi(L_1(\bu))+\psi(L_2(\bu)))_{12} \\
&= \psi(L_1(\bu))_{12} + \psi(L_2(\bu))_{12} =\omega+\omega =\omega.
\end{align*}
To prove that $\psi(\bq) \nleq \psi(\bu)$, it remains to show that $\psi(\bq)_{12} \nleq \psi(\bu)_{12}$.
For this, we shall establish that $\psi(\bq)_{12}\neq\omega$.

If $\bq=y$ for some variable $y$, then $y\notin c(L_1(\bu))$,
since $\bq\preceq\bu$ is nontrivial. By the definition of $\psi$, we have $\psi(\bq)_{12}=b$.
If $\ell(\bq)\geq3$, then by Lemma~\ref{lemma26111501}\ref{item3}, every product of three elements
in $SR_6$ equals $\infty$. Hence, all entries of $\psi(\bq)$ are $\infty$,
which means $\psi(\bq)=[\infty]_{n\times n}$. In particular, $\psi(\bq)_{12}=\infty$.

\textbf{Case 2.} $\ell(\bq)=2$.
Write $\bq=xy$. For each variable $z$, define
\[
\alpha_z=
\begin{cases}
a, &z=x,\\
c, &z\neq x,
\end{cases}
\]
and
\[
\beta_z=
\begin{cases}
\omega, &z\in c(L_1(\bu)),\\
d, &z=y\text{ and }z\notin c(L_1(\bu)),\\
b, &z\neq y\text{ and }z\notin c(L_1(\bu)).
\end{cases}
\]
Now define an assignment $\varphi$ into $\M_n(SR_6)$ by
\[
\varphi(z)=M(\alpha_z, \beta_z).
\]
Since $\alpha_z$ and $\beta_z$ are defined independently,
this assignment is well defined even when $x=y$.

We first show that $\varphi(\bu)_{12}=\omega$.
For every $z\in c(L_1(\bu))$, we have
\[
\varphi(z)_{12}=M(\alpha_z, \beta_z)_{12}=\beta_z=\omega.
\]
Next, let $st\in L_2(\bu)$. Then $s,t\notin c(L_1(\bu))$, and so
$\alpha_s\in\{a, c\}, \beta_s, \beta_t\in\{b, d\}$.
The product $\alpha_s\beta_t$ could equal $ad$ only if
$s=x$ and $t=y$. This would imply $st=xy=\bq$, contradicting
the assumption that $\bq\preceq\bu$ is nontrivial. Therefore,
\[
\alpha_s\beta_t\in\{ab, cb, cd\}=\{\omega\}.
\]
Moreover,
$\beta_sc\in\{bc, dc\}=\{\omega\}$.
Consequently,
\[
\begin{aligned}
\varphi(st)_{12}=(\varphi(s)\varphi(t))_{12}
&=\bigl(M(\alpha_s, \beta_s)M(\alpha_t, \beta_t)\bigr)_{12}\\
&\stackrel{\eqref{id26111901}}=\alpha_s\beta_t+\beta_sc
=\omega+\omega=\omega.
\end{aligned}
\]
Summing over all summands of $\bu$ yields
\[
\varphi(\bu)_{12}
=\varphi(L_1(\bu))_{12}+\varphi(L_2(\bu))_{12}=\omega+\omega
=\omega.
\]

Finally, $\alpha_x=a$, while $\beta_y=\omega$ if
$y\in c(L_1(\bu))$, and $\beta_y=d$ otherwise.
Since $a\omega=ad=\infty$, we obtain
\begin{align*}
\varphi(\bq)_{12} &= \varphi(xy)_{12} = (\varphi(x)\varphi(y))_{12}
= \bigl(M(\alpha_x,\beta_x)M(\alpha_y,\beta_y)\bigr)_{12} \\
&\stackrel{\eqref{id26111901}}= \alpha_x\beta_y+\beta_xc
= a\beta_y+\beta_xc = \infty.
\end{align*}
Since $\infty\nleq \omega$, it follows that $\varphi(\bq)_{12} \nleq \varphi(\bu)_{12}$.
Therefore, $\varphi(\bq) \nleq \varphi(\bu)$.
\end{proof}

With these preparations in place, we now provide a finite equational basis for $\M_n(SR_6)$.

\begin{proposition}\label{sr:thm:matrix}
Let $n \geq 2$ be an integer. The ai-semiring variety $\V(\M_n(SR_6))$ is defined by the identities \eqref{sr:eq:N}--\eqref{sr:eq:R}.
\end{proposition}
\begin{proof}
By Proposition~\ref{sr:lem:sound}, $\mathbf{M}_n(SR_6)$ satisfies the identities \eqref{sr:eq:N}--\eqref{sr:eq:R}.
It remains to prove that every nontrivial inequality $\bq \preceq \bu$ satisfied by $\mathbf{M}_n(SR_6)$ is derivable from \eqref{sr:eq:N}--\eqref{sr:eq:R}.

By Proposition~\ref{sr:lem:absorption}, any such inequality must satisfy at least one of the following two conditions:
either $\ell(\bu_i)\geq3$ for some $\bu_i\in \bu$, or $c(L_1(\bu))\cap c(L_2(\bu))\neq \emptyset$.
As shown in the proof of Proposition~\ref{sr:lem:absorption}, in the first case $\bq\preceq \bu$ is derivable from \eqref{sr:eq:N}, and in the second case it is derivable from \eqref{sr:eq:L} or \eqref{sr:eq:R} together with \eqref{sr:eq:N}.
Therefore, any nontrivial inequality of $\mathbf{M}_n(SR_6)$ is derivable from \eqref{sr:eq:N}--\eqref{sr:eq:R}.
\end{proof}

\begin{corollary}
Let $n\geq2$ be an integer. Then $\V(\M_n(SR_6)) = \V(\M_2(SR_6))$.
\end{corollary}
\begin{proof}
This is a direct consequence of Proposition~\ref{sr:thm:matrix}.
\end{proof}

As the final result of this section,
we determine the cardinality of the subvariety interval $[\V(SR_6), \V(\M_2(SR_6))]$,
which consists of all subvarieties of $\V(\M_2(SR_6))$ that contain $\V(SR_6)$.

For every integer $n\geq 3$, define
\[
\bs_n=\sum_{i=0}^{n-1} x_i x_{i+1}\quad {\rm and} \quad \br_n=\sum_{i=0}^{n-1} x_{i+1} x_i,
\]
where all $x_i$ are distinct variables and subscripts are taken modulo $n$.
Let $\delta_n$ denote the inequality
\[
\br_n  \preceq \bs_n.
\]

For every integer $m\geq 3$, let
$C_m$ denote the set $\{b_0, \dots, b_{m-1}\}\cup\{\omega, \infty\}$,
where the elements are pairwise distinct.
Define an addition $+$ and a multiplication $\cdot$ on $C_m$ by
\[
r+s=
\begin{cases}
r, & r=s,\\
\infty, & r\neq s,
\end{cases}
\quad
b_i\cdot b_j=
\begin{cases}
\omega, & j\equiv i+1\pmod m,\\
\infty, & j\not\equiv i+1\pmod m.
\end{cases}
\]
Additionally, every product containing $\omega$ or $\infty$ equals $\infty$.
Then $C_m$ forms an ai-semiring.
Note that $\infty$ serves as both the multiplicative zero element and the additive maximum element.
Moreover, we have
\[
b_0b_1 = b_1b_2 = \cdots = b_{m-2}b_{m-1}
= b_{m-1}b_0 = \omega.
\]
These algebras $C_m$ are examples of a class of ai-semirings introduced in \cite{gr2026}.

\begin{lemma}\label{iv:lem:Cq}
Let $m \geq 3$ be an integer. Then the ai-semiring $C_m$ lies in
$\V(\mathbf{M}_2(SR_6))$.
\end{lemma}
\begin{proof}
By Proposition~\ref{sr:thm:matrix}, it suffices to show that $C_m$ satisfies
the identities \eqref{sr:eq:N}--\eqref{sr:eq:R}.
Indeed, the product of any three elements of $C_m$ is $\infty$.
This is because the product of any two elements is either $\omega$ or $\infty$,
any product containing $\omega$ equals $\infty$, and $\infty$ is the multiplicative zero.
Since $\infty$ is the additive maximum element, it follows immediately that
$C_m$ satisfies the identity \eqref{sr:eq:N}.

Let $\varphi$ be an arbitrary assignment into $C_m$.
Since any product of three elements is $\infty$, we have $\varphi(x_1x_2x_3)=\infty$.
The product of any two elements is either $\omega$ or $\infty$, so
$\varphi(xy) \in \{\omega, \infty\}$.
If $\varphi(xy)=\omega$, then $\varphi(x)$ and $\varphi(y)$ must be $b_i$ and $b_j$
for some $0\leq i, j \leq m-1$.
Consequently,
\[
\varphi(x+xy)=\varphi(x)+\varphi(xy)=b_i+\omega=\infty=\varphi(x_1x_2x_3).
\]
If $\varphi(xy)=\infty$, then
\[
\varphi(x+xy)=\varphi(x)+\varphi(xy)=\varphi(x)+\infty=\infty=\varphi(x_1x_2x_3).
\]
Thus $C_m$ satisfies the identity \eqref{sr:eq:L}.
Similarly, $C_m$ satisfies the identity \eqref{sr:eq:R}.
\end{proof}

\begin{lemma}\label{iv:lem:cycles}
Let $m, n\geq 3$ be integers. Then $C_m$ satisfies $\delta_n$ if and only if
$m$ does not divide $n$.
\end{lemma}
\begin{proof}
Suppose first that $m$ does not divide $n$.
Let $\varphi$ be an arbitrary assignment into $C_m$.
If $\varphi(\bs_n)=\infty$, then $\varphi(\br_n) \leq \varphi(\bs_n)$, since $\infty$ is the additive maximum element.

Now suppose that $\varphi(\bs_n)\neq \infty$, then $\varphi(\bs_n)=\omega$, and so $\varphi(x_ix_{i+1})=\omega$ for all $0\leq i<n$.
We may write $\varphi(x_i)=b_{k_i}$ for $0\leq i<n$.
Then $b_{k_i}b_{k_{i+1}}=\omega$ for $0\leq i<n$, and so
$k_{i+1} \equiv k_i+1 \pmod{m}$ for $0\leq i<n$.
Note that $b_{k_n}=b_{k_0}$, since $x_n=x_0$.
We can deduce that $0 \equiv n \pmod{m}$, and so $m$ divides $n$, a contradiction.
Hence $\varphi(\bs_n)=\infty$. Consequently, $C_m$ satisfies $\delta_n$.

Conversely, assume that $m$ divides $n$.
Consider an assignment $\psi$ into $C_m$ defined by $\psi(x_i)=b_j$,
where $j$ is the remainder when $i$ is divided by $m$, $0\leq i<n$.
Since $m$ divides $n$, one can deduce that $\psi(x_ix_{i+1})=\omega$ for all $0\leq i<n$, and so $\psi(\bs_n)=\omega$.
Since $m \geq 3$, we obtain
\[
\psi(x_1x_0)=\psi(x_1)\psi(x_0)=b_1b_0=\infty,
\]
so $\psi(\br_n)=\infty$ and thus $\psi(\br_n) \nleq \psi(\bs_n)$.
Hence $C_m$ does not satisfy $\delta_n$.
\end{proof}

\begin{proposition}\label{iv:thm:sr6}
The interval $[\V(SR_6), \V(\M_2(SR_6))]$ contains continuum many varieties.
\end{proposition}
\begin{proof}
Let $\mathcal{P}$ be the set of all odd prime numbers.
Since the multiplicative reduct of $SR_6$ is commutative,
it follows that $\V(SR_6)$ satisfies the inequality $\delta_p$ for all $p \in \mathcal{P}$.
By Lemma~\ref{iv:lem:Cq}, $C_p$ lies in $\V(\M_2(SR_6))$ for all $p \in \mathcal{P}$.
For each subset $J$ of $\mathcal{P}$,
let $\mathcal{V}_J$ denote the variety generated by $SR_6$ and $C_p$, $p \in J$.
Then $\mathcal{V}_J$ contains $\V(SR_6)$ and is a subvariety of $\V(\M_2(SR_6))$.
Hence $\mathcal{V}_J$ lies in the interval $[\V(SR_6), \V(\M_2(SR_6))]$.

Let $J_1$ and $J_2$ be distinct subsets of $\mathcal{P}$.
Without loss of generality, assume that $J_1 \nsubseteq J_2$.
Then there exists $p\in J_1$ such that $p\notin J_2$.
By Lemma~\ref{iv:lem:cycles}, $\mathcal{V}_{J_2}$ satisfies $\delta_p$,
but $\mathcal{V}_{J_1}$ does not satisfy $\delta_p$.
Thus $\mathcal{V}_{J_1} \neq \mathcal{V}_{J_2}$.
We have established an injective mapping between
the power set of $\mathcal{P}$ and the interval $[\V(SR_6), \V(\M_2(SR_6))]$.
Since the power set of $\mathcal{P}$ has the cardinality of the continuum,
we conclude that $[\V(SR_6), \V(\M_2(SR_6))]$ also has the cardinality of the continuum.
\end{proof}

\section*{Acknowledgements}
The authors acknowledge the use of AI tools for exploratory discussions and language refinement during the preparation of this work. All mathematical statements, proofs, and references were independently checked, revised, and finalized by the authors, who take full responsibility for the content of this paper.

\bibliographystyle{amsplain}

\end{document}